\documentclass[11pt]{amsart}
\makeatletter
\def\@normalsize{\@setsize\normalsize{10pt}\xpt\@xpt
\abovedlayskip 10pt plus2pt minus5pt\belowdisplayskip
\abovedisplayskip \abovedisplayshortskip \z@
plus3pt\belowdisplayshortskip 6pt plus3pt
minus3pt\let\@listi\@listI}
\def\subsize{\@setsize\subsize{12pt}\xipt\@xipt}
\def\section{\@startsection {section}{1}{\z@}{1.0ex plus 1ex minus
2ex}{.2ex plus .2ex}{\large\bf}}
\def\subsection{\@startsection {subsection}{2}{\z@}{.2ex plus 1ex}
{.2ex plus .2ex}{\subsize\bf}} \makeatother

\newtheorem{lemma}{Lemma}

\newtheorem{theorem}{Theorem}

\newtheorem{corollary}{Corollary}

\theoremstyle{remark}

\begin{document}

\baselineskip=12pt

\address{Helmut-K\"autner Str. 25, 81739 Munich. Germany  }

\email{Alexei.Ostrovski@gmx.de}

\subjclass[2000]{Primary  54C10; Secondary 54H05,  54E40, 03E15.}

\keywords{  Borel sets, locally closed sets, clopen sets, open and closed functions, Borel isomorphism.}

\title{\bf {Preservation of  the  Borel class under open-$LC$  functions
}}
\author {Alexey Ostrovsky}
\maketitle

\begin{abstract}

Let $X$ be a Borel subset of the Cantor set  \textbf{C} of
additive or multiplicative   class ${\alpha},$      and $f: X \to
Y$ be a continuous  function onto $Y  \subset  \textbf{C}$  with compact preimages of points.
If the image $f(U)$ of every clopen set $U$  is the intersection  of an  open
and a closed set,   then   $Y$  is a Borel set of  the same class ${\alpha}$.
  This result
generalizes similar results for open and closed functions.

\end{abstract}

\vspace{0.15in}

\section{Introduction}

   \vspace{0.15in}

   Let $X$ be a Borel subset of the Cantor set  \textbf{C} of additive or multiplicative   class ${\alpha},$      and $f: X \to
   Y$ be
a continuous  function  onto $Y \subset  \textbf{C}$  with  compact preimages of points.

   It is well known that if the image $f(U)$ of every clopen set $U$    is an open subset  of $Y$,  then   $Y$  is a Borel set of  the same class  $\alpha$ \cite{Va}, \cite{T},   \cite{F}, \cite{JSR}.

      Analogously, if the image $f(U)$ of every clopen set $U$
 is a closed subset  of $Y$,  then   $Y $  is a Borel set of  the same class.

    The aim of  this note is to prove  (Theorem  2) that if the image $f(U)$ of every clopen  set $U$ is an intersection  of an open
and a closed set,    then   $Y$  is a Borel set of  the same class.

This fact is related to    the following problem   \cite[Problem
3.6.]{O2006}:

   \vspace{0.15in}

Find a class of open-Borel, continuous  functions that are close as   possible  
to open and closed functions and have compact preimages of points and  preserve  absolute
Borel class.

\vspace{0.15in}

 \section{  Related materials and basic definitions}

 \vspace{0.15in}

 All spaces in this paper are assumed to be metrizable and separable.

 Recall that  a subset  of a
topological space  is an  \emph{$LC$}-set or a locally closed set
if it is the intersection of an open and a closed set.

 Given an arbitrary  (not necessarily continuous)  function  $f$  we say that  it is

-\emph{open (resp. closed) } if  $f$ takes  open (resp. closed)  sets into
open (resp. closed) sets;

-\emph{open(resp. clopen)-$LC$} if $f$ takes  open (resp. clopen) sets  into $LC$-sets.

 \vspace{0.10in}

\vspace{0.05in}

\vspace{0.10in}

 \vspace{0.10in}

 The following assumptions will be  needed     throughout the paper.

 \vspace{0.10in}

We will   denote by $S_1(y)$ a  sequence with its limit point:
$$S_1(y) = \{y\} \cup \{y_i: y_i \longrightarrow y\}$$.

\vspace{0.15in}

It is easy to check that  a function  $f$   is closed
$\Leftrightarrow$    for every   $S_1(y)$,   every sequence $x_i
\in f^{-1}(y_i)$  ( $y_i   \not = y_j $ for $ i   \not = j $) has
a limit point in $f^{-1}(y)$;

 Indeed,  if   $f$   is closed   and, for  some $S_1(y)$, there is  no limit point in $f^{-1}(y)$
 for $x_i   \in f^{-1}(y_i)$, then  the image $ f(T)$ of the closed set $ T = cl_X  \{x_i  \}$
is not closed in $Y$.

 Conversely,  if  there is a closed $T    \subset  X$ for which $f(T)$ is not closed
 in $Y$, then  there is $S_1(y)$  such that $y  \not   \in  f(T)$ and  $y_i     \in  f(T)$. 
  
 Hence,  since  $T$ is  closed, every sequence of points  $x_i \in f^{-1}(y_i)  \cap T$  has no limit point
 in  $ f^{-1}(y)$.

\vspace{0.15in} Analogously, it is easy to check that a function
$f$   is open      $\Leftrightarrow$  for every   $S_1(y)$  and
every open ball $O(x)$, $x  \in f^{-1}(y)$, there are only
finitely many
 $y_i$  such that
   $f^{-1}(y_i)   \cap O(x)   = \emptyset$.

 \bigskip

  \section {Structure of  clopen-$LC$  functions  in the Cantor set   \textbf{C}   }

  \vspace{0.10in }

Let us first prove   the following theorem.

 \vspace{0.10in }

  \begin{theorem}  Let  $f: X  \to Y$ be a  clopen-$LC$  function from  a subset  $X$ of the Cantor set  $\textbf{C}$  onto $Y$ and the  inverse image of  every point $y$ be  compact.
  Then  $Y$ can be covered by countably many     subsets  $ Y_n$
  such that 
 the restrictions  $ f | f^{-1}(Y_n)$ are open functions $( n = 1,2,...)$  and the restriction  $ f | f^{-1}(Y_0)$ is a closed function.

    \end{theorem}

  Proof.
  Denote

 \vspace{0.15in }
A.  $X_n   =  \bigcup \{ f^{-1}(y) :$  there is $S_1(y)    \subset Y$ and  $ \tilde  x_y \in  \textbf{C} $ such that  there are $x_k \in  f^{-1}(y_k)$ with  $ x_k  \longrightarrow \tilde  x_y$
  and dist$( \tilde  x_y,  f^{-1}(y))  >  1/n$.

      \vspace{0.10in }
   \begin{lemma}
The restriction  $f|X_n$ is  an open function  onto $Y_n =   f(X_n)$.
 \end{lemma}

   \vspace{0.10in }

  Indeed,  suppose the opposite. Then

    \vspace{0.10in }
  B. for  some $y    \in   f(X_n) $  and  $d> 0 $,  there is $S_1(y)    \subset f(X_n) $ and $x   \in  f^{-1}(y)$
    such that   $y_k \longrightarrow y$ and  $dist(x, f^{-1}(y_k) ) > d.$

  \vspace{0.10in }

  Let us consider a countable compact set $S_2(y)$ obtained by replacing (see item   B)
   the isolated points  $y_k$  of  $S_1(y)$ by  $S_1(y_k)   \subset Y $  with  isolated points $y_{k_j}  \longrightarrow y_k$ selected according to item A.
    \vspace{0.05in }

      \vspace{0.05in }

       \vspace{0.10in }

Since $X$ lies in \textbf{C}    there is   a limit point $\tilde  x \in   \textbf{C}$ for   $\tilde  x_{y_k}$. The proof falls naturally  into two parts.

 (1)  If  $\tilde  x  \not \in X$, then  we can  take a clopen   (in   \textbf{C}) ball   $O_{\delta_1}(\tilde x)$,  $\delta_1 < 1/n$, and   a clopen   (in   \textbf{C}) ball   $O_{\delta_2}(
 x)$,
where $\delta_2 < d$,  according to B. It is clear  that   $D =
O_{\delta_1}(\tilde x)     \cup O_{\delta_2}( x)$ is a clopen set
in \textbf{C}
 and hence  $S_2(y)  \cap f(D)$  is the intersection of a closed set $F$ and an open set $U$ in $S_2(y)$.  We can suppose  that  $S_2(y)  \cap f(D)$  contains  $y,$   $y_{k_j}$ ($j,k =1,2,...),$ and obviously   $y_k \not    \in  f(D).$  Since the points  $y_{k_j}$
   are   dense in $S_2(y)$   and $y  \in U$, we obtain a contradiction
   that     $y_k     \in  f(D)$.

(2) If  $\tilde  x  \in X$, then    we can repeat   (1) for   $D =
O_{\delta_1}(\tilde x) $ in the case  $x  \in f^{-1}(y),$ else if   $x  \in X  \setminus  f^{-1}(y),$ we can repeat (1).

    \qed

        \vspace{0.10in }
   \begin{lemma}
$f$ is a closed   function at every point  of  $Y_0 = Y \setminus
\bigcup_nY_n$.   Hence, $f|X_0$ is a closed function  onto $Y_0 =
f(X_0)$.
 \end{lemma}

   \vspace{0.10in }

 Since  the   preimages of points are compact, the assertion of the lemma follows from  the definition of the sets $X_n$ in A.

 \qed

    \bigskip

      \section {Preservation of Borel class by  clopen-$LC$  functions  in the Cantor set   \textbf{C}   }

      \bigskip

\begin{theorem} Let   $f : X \to Y$  be a continuous,  clopen-$LC$  function     and the  inverse image of  any point $y$ be  compact.
  If X is a Borel set of additive or multiplicative   class ${\alpha}$  in  \textbf{ C},
  then   $Y$  is a Borel set of  the same class  in  \textbf{ C}.
 \end{theorem}
\vspace{0.15in}

To begin, we remark   that the notation below is
 as in the proof of Theorem 1.

    \vspace{0.10in }
   \begin{lemma}   The restriction
$g_n =f|f^{-1}(cl_YY_n)$ is an open function (n = 1,2,...).

 \end{lemma}

   \vspace{0.10in }

 Indeed,  suppose that      $y  \in  (cl_YY_n)  \setminus Y_n, $  so  there are  $  y_k \longrightarrow y$ with $y_k   \in Y_n.$    The method of proof of  Lemma 1  works for this case too, and a repeated application of this method as in cases   (1) and (2)  enables  us to  conclude that $g_n$ is open at every point $x  \in f^{-1} (y).$

    We  will  establish the lemma if we prove the following statement:

    Let  $   y_k \longrightarrow y$ and $y_k   \in (cl_YY_n)  \setminus Y_n. $  
    Then, for every  $x  \in f^{-1} (y)$       and every  open $O_{\delta}(x)$,  the intersection $O_{\delta}(x) \cap     f^{-1} (y_{k})$    is a nonempty set  for  infinitely many $k
    =1,2,...$.

              Suppose that this statement is false.

           Pick $y  \in Y_n$ and a clopen $V$ in $X$ that intersects $f^{ -
1} (y)$. Let $A = cl_Y Y_n$  and assume for contradiction that $f
(V ) \cap   A$  is not a neighborhood of $y$ in $A$. Since (Lemma
1), $f (V ) \cap Y_n$ is open in $Y_n,$  there are $y_n \in cl_A(f
(V )  \cap  Y_n) \setminus A$ with $y_n \longrightarrow y$, and
then $y_{n_j } \in f (V )  \cap Y_n$ with $y_{n_j} \longrightarrow
y_n.$  Thus, one gets a closed set $ F \subset   A $ with $f
(V ) \cap F $ not locally closed, which contradicts the assumption
that $ f (V )$ is locally closed.


      \qed
           \vspace{0.10in }

         Now, we turn to the proof of Theorem 2. By the above Lemma 3, we can assume that every $Y_i$ is closed in $Y$  and, hence, $Y_0$ is $G_{\delta}$ in $Y.$

           Since $f$ is continuous, every preimage $ f^{-1}(Y_n)$ $(n=1,2,...)$ is a closed subset of $X$ and  $ f^{-1}(Y_0)$ is  a  $G_{\delta}$-set of $X.$


  Hence, according to the classical  theorems on the preservation of a Borel class by closed and open  functions with compact preimages of points  \cite{Va}, \cite{T},   \cite{F}, \cite{JSR},  we find that  every $Y_n$ is an absolute Borel set of the same class.

                        If $X$  is of additive class $\alpha$, then $Y$ is of additive  class $\alpha $ because it is a countable union of the sets $Y_i.$

        Suppose  that  $X$ is a Borel set  of multiplicative   class   ${\alpha}.$ For ${\alpha} $= 1 ($X$ is  an absolute  $G_{\delta}$-set),  the conclusion follows from the recent results of S. Gao and V. Kieftenbeld \cite{SG},  P. Holicky and R. Pol \cite{HP},        and hence $Y$ is an absolute  Borel sets  of the same class         ${\alpha}.$

          Similar to  \cite[Theorem 7]{7}, we  can easily deduce  our statement for  multiplicative class          ${\alpha >1}.$

              \qed

              \vspace{0.10in }
         \section {The case of separable metric spaces.} 
         
            \vspace{0.10in }         
         
         A slight change in the proof   of Theorem 1  actually shows  the following:

          \begin{corollary}
Let  $f: X  \to Y$ be an  open-$LC$  function  such that the  inverse
image of  every point $y$ is  compact.
  Then  $Y$ can be covered by countably many     subsets  $ Y_n$
  such that 
 the restrictions  $ f | f^{-1}(Y_n)$ are open functions $( n = 1,2,...)$  and the restriction  $ f | f^{-1}(Y_0)$ is a closed function.

\end{corollary}

 The same  conclusion can be drawn from the proof of   Theorem 2:

  \begin{corollary} Let   $f : X \to Y$  be a continuous, open-$LC$ function onto $Y$ and the  inverse image of  every point $y$ be  compact.

  If X is an absolute Borel set of additive or multiplicative   class ${\alpha}$,
  then   $Y$  is an absolute  Borel set of  the same class.
 \end{corollary}
\vspace{0.15in}

A function $f$ is  a \emph{countable  homeomorphism} if  $X$  can
be partitioned into countably many pairwise disjoint sets $X_i$
such that every restriction  $f | X_i$  is a homeomorphism.

  \begin{corollary}
  Let   $f : X \to Y$   be a one-to-one function  between separable metric spaces $X$ and $Y$, such that $f$ and  $f^{-1}: Y \to X$ are
  open-$LC$ functions.
  Then   $f$ is  a countable homeomorphism.

 \end{corollary}

 Indeed, let us take the sets  $H_n^{X}  \subset X $ and
 $H_k^{Y}  \subset Y $  such that  each of $f|H_n^{X}$   and  $f^{-1}|H_k^{Y}$
 is  a one-to-one  continuous function.

 Then  every  restriction $f|X_n$, where $X_n =H_n^{X}  \cap f^{-1}(H_k^{Y})$, is a
 homeomorphism.

   \qed
            \vspace{0.15in }

If  $X$ and $Y$ are   abs. Souslin sets,  then   Corollary 3
follows from the results  of  J.E. Jayne and C.A. Rogers
\cite{JR}.

    \vspace{0.15in }


\bigskip

  \begin{thebibliography}{16}

      \bigskip

\bibitem{F} R. C. Freiwald, Images of Borel Sets and k-Analytic Sets,  {\sl Fund. Math.} 75:1 (1972), 35--46.


\bibitem{SG} S. Gao and V. Kieftenbeld,
Resovable Maps Preserve Complete Metrizability,
 {\sl Proc. Am. Math. Soc.} v.138 (2010), 2245--2252. 

\bibitem{HP} P. Holicky, R. Pol,
On a Question by Alexey Ostrovsky Concerning Preservation of
Completeness, {\sl  Topol. Appl. }157 (2010), 594--596.


\bibitem{JR} J.E. Jayne and C.A. Rogers,  First Level Borel Functions and Isomorphisms,  {\sl J. Math. Pures Appl.} 61  (1982), 177--205.



 \bibitem{7}
A.  Ostrovsky,  Finite-to-One Continuous s-Covering Mappings.  {\sl Fund. Math. }
194, (2007), 89--93.

\bibitem{O2006} A. Ostrovsky,  Maps of Borel Sets.   {\sl Proc. Steklov Inst. Math.},  v. 252  (2006), 225--247.


\bibitem{JSR} J.  Saint Raymond,
Preservation of the Borel Class under Countable-Compact-Covering
Mappings,
 {\sl Topol. Appl.},
v. 154, no.  8 (2007),  1714--1725.







\bibitem{T}
A. D. Taimanov, On Closed Mappings. IIî,  {\sl Mat. Sb. }52
(1960), 579--588. (Russian).


\bibitem{Va}
I. Vainstain,  On Closed Mappings, {\sl Mosk. Gos. Univ. Uchen.
Zap.}
{\bf 155} Mat. 5 (1952) 3--53.



 \end{thebibliography}
\end{document}